\documentclass{article}

\usepackage{arxiv}
\usepackage{amsmath}
\usepackage[utf8]{inputenc} 
\usepackage[T1]{fontenc}    
\usepackage{hyperref}       
\usepackage{url}            
\usepackage{booktabs}       
\usepackage{amsfonts}       
\usepackage{nicefrac}       
\usepackage{microtype}      
\usepackage{cleveref}       
\usepackage{lipsum}         
\usepackage{graphicx}
\usepackage{doi}

\title{Merging and disconnecting resonance tongues in a pulsing excitable microlaser with delayed optical feedback}

\author{ \href{https://orcid.org/0000-0003-3044-1846}{\includegraphics[scale=0.06]{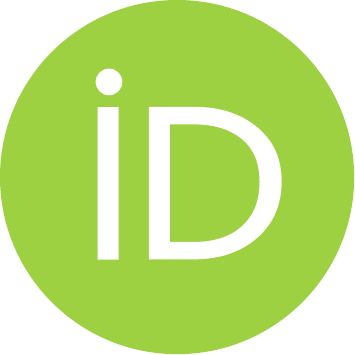}\hspace{1mm}Soizic Terrien}\thanks{soizic.terrien@univ-lemans.fr} \\
	Laboratoire d'Acoustique de l'Université du Mans (LAUM), UMR 6612,\\ Institut d'Acoustique - Graduate Schhol (IA-GS), CNRS, Le Mans Université, France\\
	\And
	\href{https://orcid.org/0000-0002-8940-230X}{\includegraphics[scale=0.06]{orcid.pdf}\hspace{1mm}Bernd Krauskopf} \\
	Department of Mathematics and Dodd-Walls Centre for Photonic and Quantum Technologies,\\ The University of Auckland,
	 Private Bag 92019 Auckland 1142, New Zealand\\
	\And
	Neil G.R. Broderick \\
	Department of Physics and Dodd-Walls Centre for Photonic and Quantum Technologies,\\
The University of Auckland, Private Bag 92019
Auckland 1142, New Zealand\\
	 \AND
	 \href{https://orcid.org/0000-0003-1304-9172}{\includegraphics[scale=0.06]{orcid.pdf}\hspace{1mm}Venkata A. Pammi} \\
	 Universit{\'e} Paris-Saclay, CNRS, Centre de Nanosciences et de Nanotechnologies,\\
 91120 Palaiseau, France \\
	 \And
	 R{\'e}my Braive \\
	Universit{\'e} Paris-Saclay, Universit{\'e} Paris Cit{\'e}, CNRS, Centre de Nanosciences et de Nanotechnologies, \\
 91120 Palaiseau, France, and Institut Universitaire de France, Paris, France \\
	 \And
	Isabelle Sagnes, Gr{\'e}goire Beaudoin, Konstantinos Pantzas and  \href{https://orcid.org/0000-0003-0092-7829}{\includegraphics[scale=0.06]{orcid.pdf}\hspace{1mm} Sylvain Barbay} \\
	Universit{\'e} Paris-Saclay, CNRS, Centre de Nanosciences et de Nanotechnologies,\\
 91120 Palaiseau, France 
}

\renewcommand{\headeright}{preprint}
\renewcommand{\undertitle}{preprint}
\renewcommand{\shorttitle}{Merging and disconnecting resonance tongues in an excitable system with delay}

\hypersetup{
pdftitle={A template for the arxiv style},
pdfsubject={q-bio.NC, q-bio.QM},
pdfauthor={David S.~Hippocampus, Elias D.~Striatum},
pdfkeywords={First keyword, Second keyword, More},
}

\begin{document}
\maketitle

%
%
%
%
%
%
%

\begin{abstract}
Excitability, encountered in numerous fields from biology to neurosciences and optics, is a general phenomenon characterized by an all-or-none response of a system to an external perturbation of a given strength. When subject to delayed feedback, excitable systems can sustain multistable pulsing regimes, which are either regular or irregular time sequences of pulses reappearing every delay time. Here, we investigate an excitable microlaser subject to delayed optical feedback and study the emergence of complex pulsing dynamics, including periodic, quasiperiodic and irregular pulsing regimes.  This work is motivated by experimental observations showing these different types of pulsing dynamics. A suitable mathematical model, written as a system of delay differential equations, is investigated through an in-depth bifurcation analysis. We demonstrate that resonance tongues play a key role in the emergence of complex dynamics, including non-equidistant periodic pulsing solutions and chaotic pulsing. The structure of resonance tongues is shown to depend very sensitively on the pump parameter. Successive saddle transitions of bounding saddle-node bifurcations constitute a merging process that results in unexpectedly large regions of locked dynamics, which subsequently disconnect from the relevant torus bifurcation curve; the existence of such unconnected regions of periodic pulsing is in excellent agreement with experimental observations. As we show, the transition to unconnected resonance regions is due to a general mechanism: the interaction of resonance tongues locally at an extremum of the rotation number on a torus bifurcation curve. We present and illustrate the two generic cases of disconnecting and of disappearing resonance tongues. Moreover, we show how a pair of a maximum and a minimum of the rotation number appears naturally when two curves of torus bifurcation undergo a saddle transition (where they connect differently).
\end{abstract}


\twocolumn
\section{Introduction}
\label{sec:level1} 

Excitability refers to the spiked or pulsed response of a system at rest to an external perturbation when the perturbation amplitude exceeds the so-called excitable threshold, while no response occurs for smaller perturbations \cite{izhikevich2000neural}. This general phenomenon has been described in a variety of scientific fields from biology \cite{WedgwoodJRS21} to neurosciences \cite{IzhikevichBook} and optics \cite{turconi2013control}, and typically results from the interplay between different internal timescales in the system. After a (short) excitable pulse has been triggered, the system enters a refractory period, during which it is either impossible or much more difficult to trigger another excitable response\cite{SelmiPRL14}; note that the refractory period is significantly larger than the duration of the pulse itself. In the presence of delayed feedback, an excitable system can regenerate its own output: after the first excited pulse, the output is reinjected after the feedback delay time $\tau$ which triggers the next pulse. As the process repeats, this results in a periodic pulsing regime whose period is directly related and close to the delay $\tau$. This general mechanism for self-pulsations has been demonstrated in a variety of optical systems \cite{GarbinNC15,RomeiraNSR16,TerrienPRA17}, as well as in an excitable biological cell \cite{WedgwoodJRS21}. 

From a more general point of view, the introduction of a delay to an excitable system can induce a wealth of complex dynamics beyond regular self-pulsing \cite{garbin2017interactions,TerrienSIADS17,seidel2022influence}. This includes a high degree of multistability \cite{yanchuk2009delay} and an enhanced dynamical complexity, such as, for example, quasiperiodic \cite{munsberg2020topological} and chaotic regimes. In particular, recent experimental and numerical investigations have demonstrated that an excitable microlaser with delayed feedback can sustain multistable periodic pulsing regimes \cite{TerrienOL18,TerrienPRR20}. Depending on the ratio between the internal timescales of the excitable microlaser and on the delay time of the feedback loop or external cavity, these include pulsing patterns with equidistant pulses of equal amplitude or with non-equidistant pulses of different amplitudes \cite{TerrienPRE21}. The emergence of equidistant pulsing patterns is well understood \cite{yanchuk2009delay,KrauskopfWalker,TerrienPRR20}. The emergence of non-equidistant pulsing regimes, on the other hand, has been suggested to originate in resonance phenomena associated with locked periodic orbits on stable tori \cite{TerrienPRE21} --- yet this was still to be investigated, which motivated the work presented here.

We adopt a dynamical systems point of view to investigate experimentally and numerically the emergence of multi-frequency dynamics in an optical realisation of an excitable system with delayed feedback. We consider an excitable microlaser subject to delayed optical feedback \cite{SelmiPRL14,ElsassEPJD10}, whose study is motivated by potential applications to neuromorphic photonic computing \cite{pammi2019photonic, shastri2021photonics}. Its feedback-induced dynamics is investigated both experimentally and with a suitable mathematical model --- the Yamada equations with delayed feebback \cite{Yamada93,KrauskopfWalker}, which take the form of three coupled delay differential equations (DDEs) with two slow and one fast variables. Compared to ordinary differential equations, solving DDEs requires specific numerical methods due to their infinite-dimensional nature \cite{roose2007continuation}. We use the Matlab-based numerical continuation software DDE-Biftool \cite{Engelborghs00report,EngelborghsACM02,sieber2014dde} to perform an in-depth bifurcation analysis in three parameters of practical importance: the feedback delay $\tau$, the feedback strength $\kappa$ and the pump parameter $A$. The results presented here unveil how very large, experimentally accessible locking regions emerge in the parameter space, leading to an increased and observable degree of multistability. We show that this phenomenon involves several transitions in the $(\tau,\kappa)$-plane that change the structure of regions of locked dynamics. These occur in such a very small $A$-interval that the associated switch in observed pulsing may seem instantaneous from an experimental perspective. From the mathematical point of view, on the other hand, we are able to distinguish and identify these transitions. In particular, we find two generic cases of resonance tongues interacting locally at an extremum of the rotation number on a torus bifurcation curve. These lead to the disconnection and disappearance of resonance tongues in a parameter plane, respectively, as a third parameter is changed. As we also show, extrema of the rotation number emerge naturally, including the Yamada model with delayed feedback, when two torus bifurcation curves reconnect differently at a saddle transition.

The article is organised as follows. The experimental device is described in Section~\ref{sec:device}, and experimental observations are presented and discussed. Background on the mathematical model is provided in Section~\ref{sec:mdl}. In Section~\ref{sec:2}, multi-frequency dynamics, including non-equidistant periodic pulsing regimes, quasiperiodic regimes and chaos are investigated through time-domain simulations of the model. In Section~\ref{sec:2D}, a bifurcation analysis demonstrates that resonance tongues in the $(\tau,\kappa)$-plane play a key role in the emergence of the observed multi-frequency dynamics. The sensitivity of the structure of resonance tongues to the experimentally relevant pump parameter $A$ is investigated in detail in Section~\ref{sec:pump}, and the phenomenon of disconnecting and disappearing resonance tongues is the focus of Section~\ref{sec:merge_dis}. We draw some conclusions in Section~\ref{sec:concl}.

\subsection{Experimental device and observed pulsing regimes}
 \label{sec:device}

In our optical realisation we use a  semiconductor micropillar laser with integrated saturable absorber, consisting of a 5$\,\mu$m diameter pillar laser with an original design \cite{SelmiPRL14,ElsassEPJD10,BarbayOL11}. In particular, its microcavity includes both a gain and a saturable absorber section. This microlaser emits light perpendicularly to its surface at the cavity resonance wavelength of 980\,nm and is optically pumped at around 800\,nm. The pump intensity is set just below the self-pulsing threshold such that, in the absence of feedback, the microlaser is in the excitable regime \cite{SelmiPRL14, SelmiPRE16,ErneuxPRE18}. Short optical pulse perturbations (of 80\,ps duration) can be sent by an external Ti:Sa mode-locked laser to trigger excitable responses that consist of optical pulses of approximately 200\,ps duration. To realise an external cavity providing the delayed feedback, part of the signal emitted by the microlaser is transmitted through a beamsplitter (R/T=70/30) and reflected back by a distant mirror after passing through a 5\,cm focal length lens. The resulting optical feedback delay time $\tau$ can be set to between $\sim$ 5 to 10\,ns. The light reflected by the beamsplitter is analyzed using an avalanche photodiode, amplified with a high bandwidth RF-amplifier and recorded with an oscilloscope. The microlaser is temperature controlled close to room temperature thanks to a Peltier cooler.  

In the presence of feedback, it has been shown that a first excitable pulse can regenerate itself when reinjected in the micropillar after the delay $\tau$, provided that the feedback strength is sufficiently large \cite{TerrienPRR20}. A periodic pulsing regime results with fundamental repetition frequency close to $\tau^{-1}$; this regime can coexist with harmonic pulsing solutions with several regularly timed pulses in the feedback cavity. It has been shown that individual pulse trains can be added or erased by single external optical perturbations \cite{TerrienOL18,TerrienPRR20}.  Under certain experimental conditions, non-regular pulse trains can also be emitted following a seemingly pulse-timing symmetry-breaking phenomenon \cite{TerrienPRE21}.	

\begin{figure}[t]
\includegraphics[width=\linewidth]{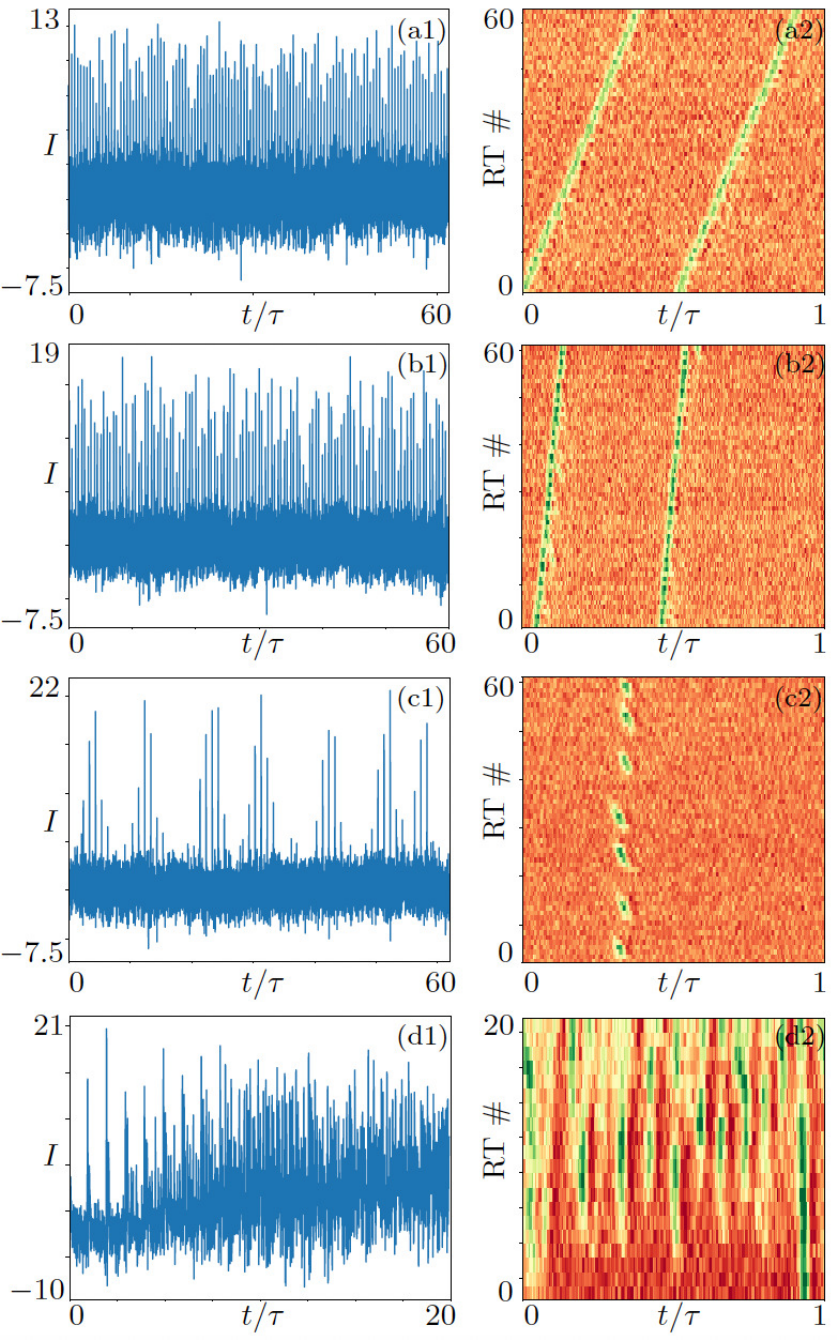}
\caption{Experimental intensity time series (left, a1-d1) and their pseudo-space representation (right, a2-d2) for: (a) equally spaced two pulses per feedback roundtrip ($\tau$ = 5.47 ns), (b) symmetry-broken two pulses per roundtrip ($\tau$ = 8.22 ns), (c) modulated quasiperiodic regime ($\tau$ = 4.78 ns), and (d) complex dynamics above the lasing threshold
($\tau$ = 4.77 ns).}
\label{fig:expe}
\end{figure}

Figure~\ref{fig:expe} illustrates the diversity of pulsing dynamics observed experimentally. The left column represents time series of the measured intensity $I$ and the right column shows their pseudo-space representation \cite{ArecchiPRA92}: temporal traces are folded at (approximately) multiples of the delay time $\tau$ and stacked vertically. In this representation, the $x$-axis represents the delay line (realized by the external feedback cavity) and the $y$-axis the recorded number of roundtrips in the external feedback cavity.  Figure~\ref{fig:expe}(a) shows a periodic regime with two equidistant pulses per feedback roundtrip, triggered by an appropriate sequence of two successive external perturbations \cite{TerrienPRR20}. The amplitude of the pulses is quite irregular, while the interpulse timing repeats consistently roundtrip after roundtrip. Such laser intensity fluctuations can be attributed to pump noise and also to detection noise since the microlaser has a very low output power \cite{TerrienPRA17}. In Figure~\ref{fig:expe}(b), the system is in the symmetry-broken regime \cite{TerrienPRE21}: it sustains two pulses per roundtrip with unequal but well-defined interpulse timings. In the pseudo-space representation, it appears as two non-equidistant pulses in the external feedback cavity.  Figure~\ref{fig:expe}(c) illustrates a regime reminiscent of quasiperiodic dynamics. In particular, it displays a strong modulation of the pulse amplitude. Note that the pulse timing is also affected, as the pseudo-space representation shows: the group of pulses shows a negative shift (towards the left) when the pulse amplitude increases. This is explained by the strong amplitude-timing coupling in excitable systems: the response time to an external perturbation gets shorter when the perturbation amplitude is increased \cite{TerrienOL18}. Importantly, all the regimes shown in panels (a-c) are multistable and coexist with the off-state (non-lasing equilibrium) of the laser. Finally, Figure~\ref{fig:expe}(d) illustrates a complex pulsing regime reminiscent of chaotic dynamics. It is recorded for a value of the pump parameter slightly above the first lasing threshold at which the off-state of the laser loses stability. As such, this regime is not triggered by external pulse perturbations, but rather by noise. The pseudo-space representation clearly shows multiple competing pulses in the feedback cavity. This complex regime cannot be observed for a long period of time because the microlaser heats up, which results in it switching off. 

Overall, Figure~\ref{fig:expe} clearly demonstrates the existence of different periodic, quasiperiodic and more complex regimes of the excitable micropillar laser with delayed optical feedback. It should be noted that the different regimes in Figure~\ref{fig:expe} are observed for different experimental parameters, in particular, for different values of feedback delay and pump intensity, which are easily tuned experimental parameters. Moreover, different microlasers (of the same design and with the same specifications) were considered, which possibly display slightly different internal parameters (including different recombination rates of carriers in the gain and absorber sections) due to fabrication.

\subsection{Background on the Yamada rate equations with delay \label{sec:mdl}}

The dynamics of the experimental system under consideration is investigated by using a straightforward extension of the original Yamada rate equation model \cite{Yamada93}, a well known system of three ordinary differential equations (ODEs) for single-mode, Q-switched lasers. It has been studied extensively, in particular, through a complete numerical bifurcation analysis which highlighted all its possible dynamics \cite{DubbeldamOC99,otupiri2020yamada}. This showed that this model exhibits an excitable regime for a large range of parameters, below the lasing threshold at which the non-lasing (\textit{i.e.}, zero-intensity) equilibrium becomes unstable.  We consider here the Yamada model with an additional delayed optical feedback term \cite{KrauskopfWalker}.  This model has been shown to produce a wealth of pulsing regimes \cite{KrauskopfWalker,TerrienSIADS17} and to describe accurately a range of dynamics observed experimentally in the micropillar laser with integrated saturable absorber and delayed optical feedback considered in this article \cite{TerrienPRA17,TerrienOL18, TerrienPRR20,TerrienPRE21}. This includes a variety of stable pulsing periodic regimes with different numbers of equidistant and non-equidistant pulses in the feedback cavity. Overall, the model is written as a system of three coupled DDEs for the dimensionless gain G, absorption Q and intensity I:
\begin{equation}
\begin{split}
\dot G &= \gamma_G(A-G-GI);\\
\dot Q &=\frac{\gamma_G}{\sigma}  (B-Q-a\sigma QI);\\
\dot I &= (G-Q-1)I+\kappa I(t-\tau).
\end{split}
\label{eq:yam}
\end{equation} 
Here, A is the pump parameter, B describes the linear absorption, a is the saturation parameter, $\gamma_G$ is the recombination rate of carriers in the gain section, and $\sigma$ is the ratio between the recombination rate of carriers in the gain and absorber sections; furthermore, time is rescaled in \eqref{eq:yam} to the photon lifetime in the cavity (which is on the order of 1-2 ps). Importantly, $\gamma_G$ is usually small in semiconductor lasers and $\sigma$ is between 0.5 and 2 for the microlaser we consider. As such, system \eqref{eq:yam} is a slow-fast dynamical system with two slow variables $G$ and $Q$, and one fast variable $I$. In the intensity equation, the delayed term describes the incoherent delayed optical feedback with feedback strength $\kappa$ and delay time $\tau$.  The influence of the feedback parameters $\kappa$ and $\tau$ has been investigated through an extensive bifurcation analysis \cite{KrauskopfWalker,TerrienSIADS17}. This highlighted, in particular, an important and increasing level of multistability as the delay $\tau$ is increased, with a large number of coexisting stable periodic pulsing solutions. 

The parameters describing material properties of the laser are fixed here throughout to $B= 2$, $\gamma_G=0.01$, $a=5.5$, and $\sigma=1.8$; these values are chosen to match the parameters considered in previous work \cite{TerrienPRE21}. The pump parameter $A$, the feedback strength $\kappa$ and the feedback delay $\tau$, on the other hand, can be changed during the experiment and are considered as bifurcation parameters. Importantly, for all the parameters combinations considered in this article, the solitary laser for $\kappa=0$ (\textit{i.e.}, the model without the feedback term) is in the excitable regime: the non-lasing equilibrium, corresponding to the laser off-state, is stable but the system can release a single intensity pulse when subject to an external perturbation with sufficiently large amplitude \cite{DubbeldamOC99}.

\section{Multi-frequency dynamics in time-domain simulations}
\label{sec:2}

\begin{figure}
\includegraphics[width=\linewidth]{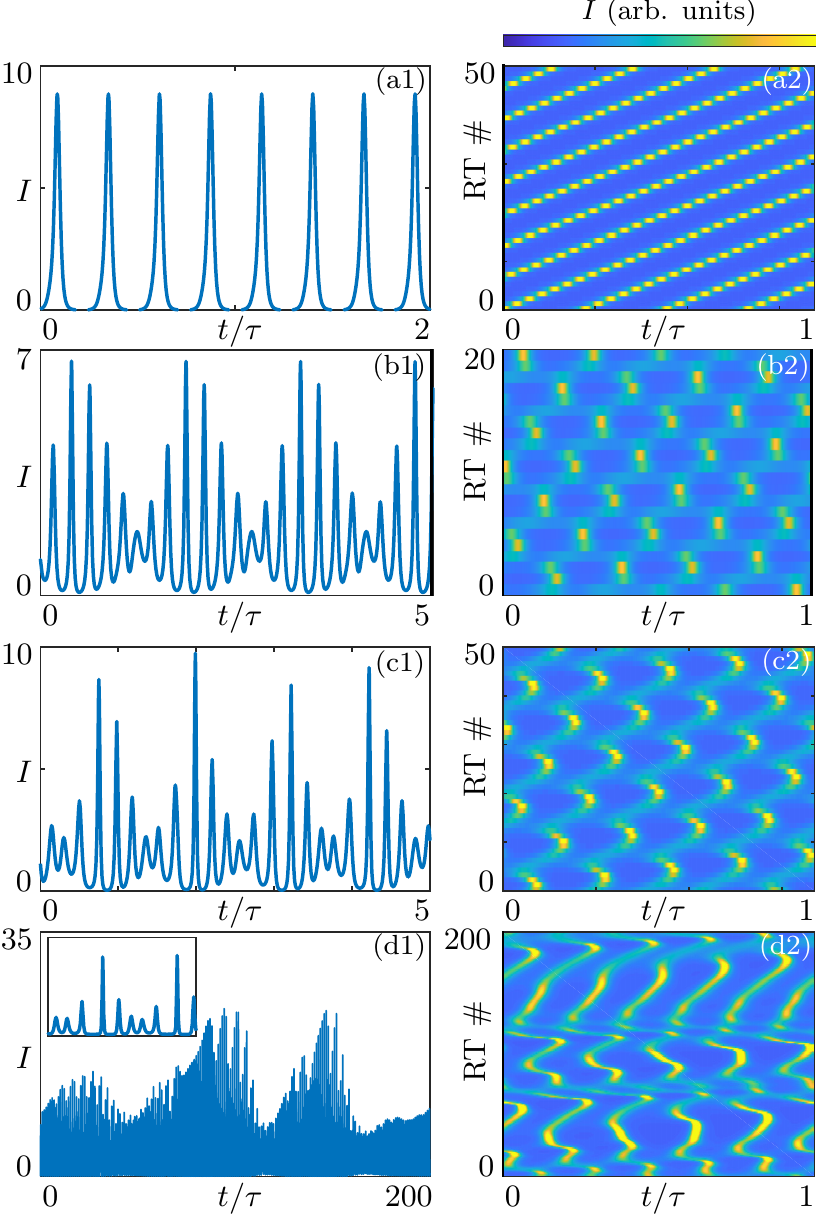}
\caption{Time-domain simulations of \eqref{eq:yam} for $A=2.7$ and $\tau=335$, showing intensity time series (left, a1-d1) and their pseudo-space representation (right, a2-d2); here $\kappa=0.007$ in (a), $\kappa=0.0159$ in (b), $\kappa=0.04$ in (c), and $\kappa=0.1$ in (d). The inset in panel (d1) shows an enlargement of the intensity time series over two roundtrips in the external feedback cavity.} 
\label{fig:TS_simu}
\end{figure}

The experimentally observed regimes in Figure~\ref{fig:expe} can be identified in time-domain simulations of the Yamada model with feedback \eqref{eq:yam}. Figure~\ref{fig:TS_simu} illustrates such simulations with time series of the intensity $I$ for fixed feedback delay $\tau=335$ and pump parameter $A=2.7$, and the four values of the feedback strength $\kappa=0.007$, $\kappa=0.0159$, $\kappa=0.04$ and $\kappa=0.1$. We checked that, for the considered value of $A$, the solitary laser (for $\kappa=0$) is indeed in the excitable regime. Moreover, in the presence of feedback, the non-lasing (\textit{i.e.}, zero-intensity) equilibrium solution of \eqref{eq:yam} is still stable for all the values of $\kappa$ considered in Figure~\ref{fig:TS_simu}: hence, the intensity $I$ remains zero in the absence of external perturbations. In the simulations, an initial external perturbation is accounted for by setting initial conditions (given for a DDE by a history segment over $[-\tau,0]$) with suitable non-zero intensity. Specifically, these are set to the (unstable) equilibrium solution of \eqref{eq:yam} with a non-zero intensity $I$, which corresponds to the continuous-wave regime of the laser \cite{TerrienSIADS17}.

The time series in Figure~\ref{fig:TS_simu}(a1-d1) are shown after a few dozens or hundreds of roundtrips to ensure that any transient dynamics has died down. The displayed dynamical regimes thus correspond to stable pulsing regimes of \eqref{eq:yam}. Panels~\ref{fig:TS_simu}(a2-d2) show the corresponding pseudo-space representation of these time series\cite{ArecchiPRA92}, as explained above. These simulation results illustrate the diversity of stable pulsing regimes observed over a small range of the feedback strength $\kappa$. 

Figure~\ref{fig:TS_simu}(a) for $\kappa=0.007$ shows a periodic pulsing regime with four equidistant pulses in the external feedback cavity, that is, over the span of one delay time $\tau$ that constitutes the feedback loop. Note that the amplitude also repeats exactly here because \eqref{eq:yam} does not feature (pump or other) noise. The positive slope of the intensity pulse trains in the pseudo-space representation in panel (a2) shows that the period of pulsing is slightly larger than $\tau/4$, which is due to the latency time of the system to the re-injected perturbation  \cite{yanchuk2009delay,TerrienSIADS17,KrauskopfWalker}.

Figure~\ref{fig:TS_simu}(b-c) for increasing values of $\kappa$ (where all the other parameters are fixed) show examples of dynamics on a torus, which may be quasiperiodic or locked to an attracting periodic orbit. Locked periodic solutions are found inside resonance regions or \emph{resonance tongues}, which are regions of a parameter plane that are bounded by curves of saddle-node bifurcations of periodic orbits. Resonance tongues emerge from resonance points, which are points along a torus (or Neimark-Sacker) bifurcation curve where the rotation number is rational \cite{Kuznetsov2013elements,keane2015delayed}. Quasiperiodic solutions, on the other hand, are found in the parameter plane along curves `in between' infinitely many and generally very narrow resonance tongues. Figure~\ref{fig:TS_simu}(b) for $\kappa=0.0159$ is past a torus bifurcation and shows a locked periodic orbit on a torus. Note that the period of this periodic regime is not directly related to an integer submultiple of the feedback delay time $\tau$. Rather, a single period of the periodic solution displays several intensity pulses with different amplitudes, which is characteristic of a locked periodic orbit on a stable invariant torus; see panel (b1). Figure~\ref{fig:TS_simu}(c) for $\kappa= 0.04$ shows a quasiperiodic regime (or a periodic regime with a very large period); in particular, the time series does now not repeat exactly but displays roughly five pulses in the feedback loop, with a deeply modulated pulse amplitude. The pseudo-space representation in panel (c2) highlights the strong amplitude-time coupling of the pulses \cite{TerrienOL18}. The fact that a stable quasiperiodic regime is observed in panels (c) for a nearby value, yet larger value of $\kappa$ support the interpretation that panels (b) indeed show a locked periodic solution on a stable torus.

Finally, Figure~\ref{fig:TS_simu}(d) for the yet larger value of $\kappa= 0.1$ shows an example of a chaotic regime. Here no clear structure is observed in neither the repetition rate nor the pulse amplitudes, which is visually clear especially in panel (d2); see also panel (d1) and the enlargement of the intensity time series over two roundtrips in the inset; note, in particular, that the modulation of the amplitude is no longer periodic.

\begin{figure}
\includegraphics[width=\linewidth]{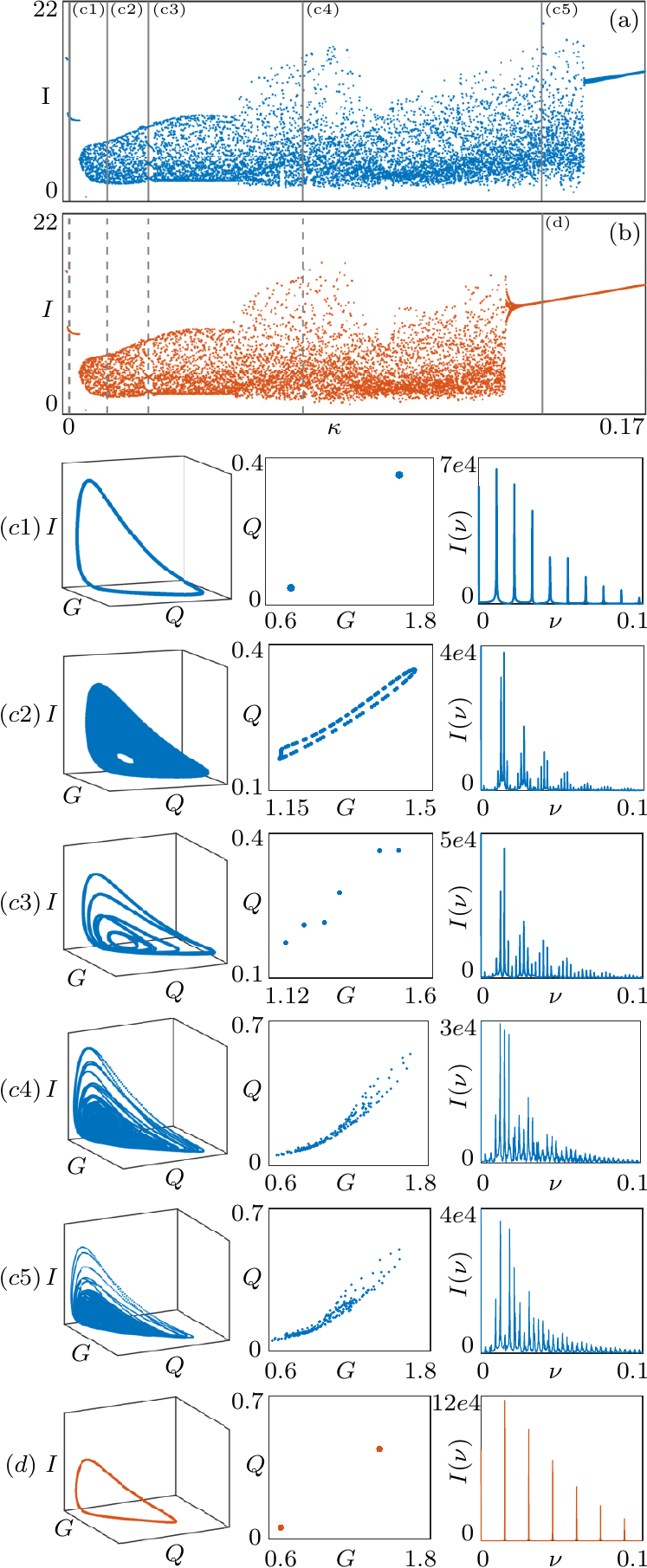}
\caption{One-parameter bifurcation diagrams of \eqref{eq:yam} for $A=2.7$ and $\tau=335$ where $\kappa$ is increased (a) and decreased (b). Also shown are representative trajectories in $(G,Q,I)$ space (left), their intersection with a Poincar{\'e} section represented by the $(G,Q)$-plane (middle), and corresponding RF spectra (right), for $\kappa=0.002$ (c1), $\kappa=0.013$ (c2), $\kappa=0.025$ (c3), $\kappa=0.07$ (c4) and $\kappa=0.14$ (c5,d).}
\label{fig:diag1D}
\end{figure}

Transitions between the different stable pulsing regimes in Figure~\ref{fig:TS_simu} as the feedback strength $\kappa$ is changed are explored further in Figure~\ref{fig:diag1D}. Panels (a) and (b) show one-parameter bifurcation diagrams for increasing and decreasing values of $\kappa$, respectively, where observed dynamics is represented by  the pulses amplitudes observed over one single roundtrip of duration $\tau$. These diagrams are obtained from sweeped simulation as follows. For the smallest values of $\kappa$, the initial condition/history is set to the (unstable) lasing equilibrium solution; for each subsequent value of $\kappa$, the previously calculated solution (\textit{i.e.}, for a slightly smaller or larger value of $\kappa$ in Figure~\ref{fig:diag1D}(a) and (b), respectively) is considered as initial history. In all cases, the simulation is run over several hundreds of roundtrips so that any transient phenomenon are disregarded, before the pulses amplitudes over one single roundtrip are recorded. Figure~\ref{fig:diag1D} also illustrates individual dynamics at selected values of $\kappa$ that are indicated by vertical lines in panels (a) and (b). In each case, we show the attractor in $(G,Q,I)$ space and in a Poincar{\'e} section represented by the $(G,Q)$-plane. It is important to note here that, due to the infinite-dimensional nature of DDEs, these are both projections: onto the three-dimensional physical space of the variables $G$, $Q$ and $I$ for the trajectory and, similarly, onto the $(G,Q)$-plane for the Poincar{\'e} map, which is defined by the fixed value $I_{eq}$ of the intensity at the unstable equilibrium point \cite{KrauskopfWalker}; specifically, $I_{eq}=1.383$ in panel (c1), $I_{eq}=1.412$ in panel (c2), $I_{eq}=1.444$ in panel (c3), $I_{eq}=1.572$ in panel (c4), and $I_{eq}=1.796$ in panels (c5) and (d). We also show the corresponding radio-frequency (RF) spectrum (the power spectral density) of the real-valued intensity time series. 

As $\kappa$ is varied, several transitions are observed between the pulsing regimes displayed in Figure~\ref{fig:TS_simu}. For the smallest values of $\kappa$ in Figure~\ref{fig:diag1D}(a), a single value is observed for the pulse amplitudes. This corresponds to a periodic regime with a fixed number of equidistant pulses of equal amplitude per roundtrip \cite{TerrienPRE21}. This periodic regime is illustrated in panel (c1) for $\kappa=0.002$, where one observes a periodic orbit in $(G,Q,I)$ space and two corresponding points in the $(G,Q)$-plane representing the Poincar{\'e} section; the RF spectrum is discrete, featuring a main peak with large harmonics owing to the strongly pulse-like nature of the oscillation. Figure~\ref{fig:diag1D}(a) shows that increasing $\kappa$ leads to a transition from this periodic to a regime with an increasingly stronger modulation of the pulse amplitude. This is clear evidence of a torus bifurcation with subsequent quasiperiodic or locked dynamics. Panels (c2) for $\kappa = 0.013$ show the corresponding stable torus in projection onto $(G,Q,I)$-space; its intersection with the Poincar{\'e} section is a densely filling closed curve in the $(G,Q)$-plane, and its RF spectrum displays a large number of peaks at incommensurate frequencies; compare with Figure~\ref{fig:TS_simu}(b). Hence, the dynamics on the torus is quasiperiodic (or of very high period). When $\kappa$ is increased in Figure~\ref{fig:diag1D}(a), different (usually very narrow) resonance tongues are briefly crossed, where the dynamics is locked to a periodic orbit on the stable torus; see also Figure~\ref{fig:TS_simu}(c). Figure~\ref{fig:diag1D}(c3) shows such a locked periodic regime for $\kappa=0.025$, which lies in a larger resonance tongue corresponding to a wider range of $\kappa$. The trajectory in $(G,Q,I)$-space shows a closed trajectory, which does not cover the entire stable torus but rather winds six times around it. The Poincar{\'e} section in the $(G,Q)$-plane shows the corresponding six intersection points, which confirms that the system is in a 1:6 locking region. When $\kappa$ is increased further, a transition is observed in Figure~\ref{fig:diag1D}(a) from quasiperiodic (or high-period) dynamics on a torus to a chaotic regime; this is detected by a sudden growth of the maximal recorded amplitude. The chaotic regime is illustrated further in Figure~\ref{fig:diag1D}(c4-c5); compare with the chaotic intensity time series in Figure~\ref{fig:TS_simu}(d). Here, neither the (projected) trajectory nor the Poincar{\'e} section display a clear structure; moreover, the RF spectrum is characteristic of a chaotic regime: although some clear peaks are visible, there are multiple frequency bands in which no clear structure can be seen. This chaotic regime is observed over a large range of $\kappa$ in Figure~\ref{fig:diag1D}(a) before the system jumps back to a stable periodic pulsing regime as $\kappa$ is increased further. 

Starting from this periodic solution, Figure~\ref{fig:diag1D}(b) shows the one-parameter bifurcation diagram of \eqref{eq:yam} for decreasing $\kappa$. The comparison between the bifurcation diagrams in panels (a) and (b) highlights a region of multistability, where the chaotic regime displayed in Figure~\ref{fig:diag1D}(c5) coexists with the stable periodic regime shown in panels (d). As can be seen from its RF spectrum, this periodic regime features a different number of equidistant pulses in the feedback cavity compared to the periodic regime in panel~(c1). It is important to note that the non-lasing equilibrium solution of \eqref{eq:yam} is stable as well over the entire range of $\kappa$ considered in Figure~\ref{fig:diag1D}(a,b). Which regime is observed in practice in the presence of the overall multistability depends on the considered initial condition/history \cite{TerrienPRR20}.

\section{Bifurcation analysis in the ($\tau,\kappa$)-plane}
\label{sec:2D} 

The complex dynamics highlighted in Figures~\ref{fig:TS_simu} and~\ref{fig:diag1D} is now investigated further by means of numerical bifurcation analysis, where the feedback delay $\tau$ and the feedback strength $\kappa$ are considered as the main two bifurcation parameters. This is performed with the Matlab-based continuation software DDE-Biftool \cite{EngelborghsACM02,Engelborghs00report,sieber2014dde}, which allows for the continuation of families of equilibrium and periodic solutions of delay-differential equations, as well as their bifurcations. The results are presented in the form of bifurcation diagrams consisting of bifurcation sets in the $(\tau,\kappa)$-plane, formed by curves of codimension-one bifurcations that bound regions of different dynamics.

\subsection{Bifurcation diagram for $A=2.7$}
\label{sec:bifdiagA2d7}

\begin{figure}
\begin{center}
\includegraphics[width=0.99\linewidth]{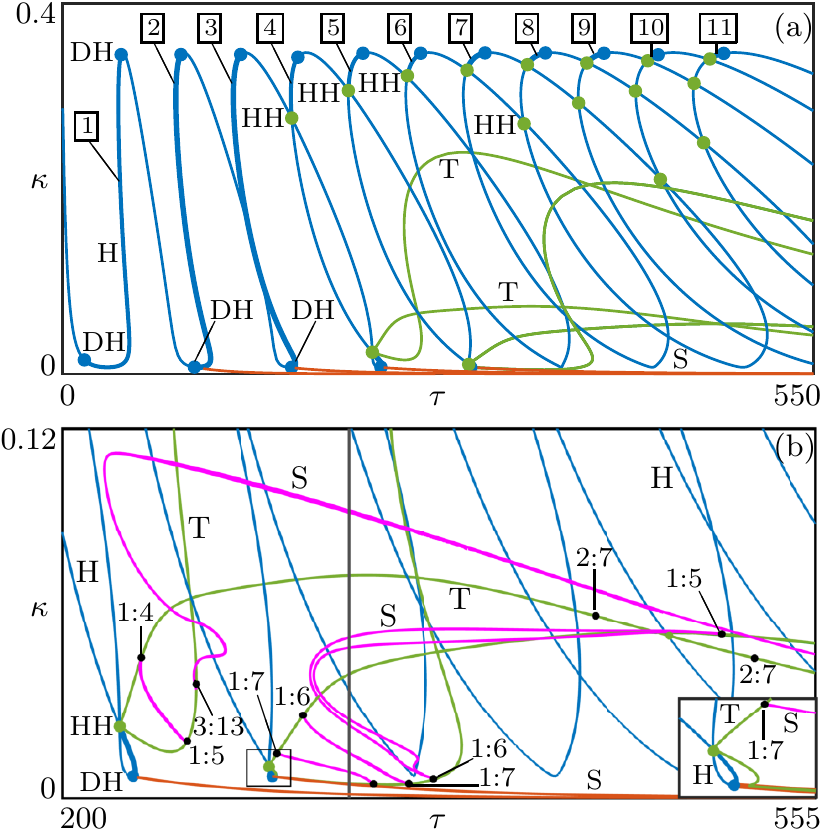}
\caption{Two-parameter bifurcation diagram of \eqref{eq:yam} in the ($\tau$,$\kappa$)-plane for $A=2.7$ (a), showing a curve H of Hopf bifurcation (blue), curves T of torus bifurcation (green), a curve S of saddle-node bifurcations of periodic orbits (red), as well as codimension-two points DH of degenerate Hopf bifurcation points (blue dots) and HH of Hopf-Hopf bifurcation (green dots). The enlargement in panel (b) additionally shows some of the (lowest-order) resonance points along the torus bifurcation curves and the curves S of saddle-node bifurcations of periodic orbits (magenta) bounding the corresponding resonance tongues. The inset shows a further enlargement of the framed area in panel (b), highlighting the emergence of the torus bifurcation curve from a Hopf-Hopf bifurcation point.}
\label{fig:2D_A2d7}
\end{center}
\end{figure}

Figure~\ref{fig:2D_A2d7} shows the bifurcation diagram in the ($\tau$,$\kappa$)-plane for the value $A=2.7$ of the pump parameter considered in Figures~\ref{fig:TS_simu} and~\ref{fig:diag1D}; its panel (a) shows the range of $\tau$ from $0$ to $550$, and panel (b) is an enlargement that focuses on the resonance structures we consider here. A main feature of Figure~\ref{fig:2D_A2d7} is a single curve H of Hopf bifurcation of the lasing (\textit{i.e.}, non-zero intensity) equilibrium of \eqref{eq:yam}: when this curve is crossed, for instance increasing $\tau$ at constant $\kappa$, a small-amplitude periodic solution emerges. More precisely, the bifurcating periodic solution is locally stable when the Hopf bifurcation is supercritical, indicated by bold parts of the curve H, and it is locally unstable along thin parts of H where the Hopf bifurcation is subcritical. The criticality  of H changes at codimension-two degenerate Hopf bifurcation points DH, indicated by blue dots in Figure~\ref{fig:2D_A2d7}. We also find codimension-two Hopf-Hopf bifurcation points HH in Figure~\ref{fig:2D_A2d7}; they arise here due to self-intersections of the curve H and also lead to changes of criticality of the Hopf bifurcation. When $\tau$ is increased for a (sufficiently large) fixed value of $\kappa$, several stable periodic solutions emerge successively from the supercritical parts of the curve H, and their periods are close to submultiples of the delay $\tau$ \cite{KrauskopfWalker,TerrienSIADS17}. Past the Hopf bifurcations where they emerge, these periodic solutions correspond to pulsing patterns with different numbers of equidistant pulses in the feedback loop \cite{TerrienSIADS17}; these numbers are shown in boxes in Figure~\ref{fig:2D_A2d7}(a) for the respective supercritical parts of the curve H. The increasing level of multistability with the delay $\tau$ is typical for delay systems \cite{yanchuk2009delay}, and has been discussed in the literature for this particular system \cite{ruschel2020limits}.

The codimension-two points DH and HH in Figure~\ref{fig:2D_A2d7} also give rise to additional bifurcation curves in the ($\tau$,$\kappa$)-plane. From each degenerate Hopf bifurcation point DH, where the Lyapunov coefficient (third-order normal form coefficient) vanishes, emerges a curve S of saddle-node bifurcation of periodic orbits \cite{KrauskopfWalker,TerrienSIADS17}; for sake of keeping this exposition focused, Figure~\ref{fig:2D_A2d7} only shows the curve S that emerges from a point DH for low $\kappa$ near $\tau = 300$, which we require for the discussion in Section~\ref{sec:pump}. Moreover, at each Hopf-Hopf bifurcation point HH, where the linearization of system \eqref{eq:yam} at the bifurcating lasing equilibrium has two pairs of complex conjugate eigenvalues on the imaginary axis, two curves of torus bifurcation T can typically emerge \cite{Kuznetsov2013elements}. For sake of clarity, only two pairs of curves T are shown in Figure~\ref{fig:2D_A2d7}: those emerging from the two points HH for low values of $\kappa$ around $\tau=335$. Indeed, many Hopf-Hopf bifurcation points are encountered as the delay $\tau$ is increased, resulting in a very complex bifurcation diagram with many more curves of torus bifurcation; this has been discussed in the literature \cite{KrauskopfWalker,TerrienSIADS17,ruschel2020limits} and is beyond the scope of this article.

\subsubsection{Connecting resonance tongues}
\label{sec:connecting}

We now focus on the two shown pairs of torus bifurcation curves T.  They are shown in Figure~\ref{fig:2D_A2d7}(b) in the enlarged parameter area of interest, where they are clearly seen to emerge from the two respective Hopf-Hopf bifurcation points; see also the further enlargement in the inset. From each torus bifurcation curve T bifurcates a smooth invariant torus on which the (multi-frequency) dynamics is either quasiperiodic or locked. Resonance tongues in the ($\tau$,$\kappa$)-plane emerge from $p\!\!:\!\!q$ resonance points, some of which are highlighted on the curves T by black dots: at these points, the pair of critical complex conjugate Floquet multipliers on the unit circle $e^{\pm 2\pi\alpha}$ has the rational rotation number $\alpha=\frac{p}{q}$, where $p$ and $q$ are relatively prime \cite{Kuznetsov2013elements}. A pair of curves S of saddle-node bifurcation of periodic orbits emerge from each $p\!\!:\!\!q$ resonance point; they form the boundaries of a $p\!\!:\!\!q$ resonance tongue in which the dynamics is $p\!\!:\!\!q$ locked on the torus, that is, periodic rather than quasiperiodic; see Figure~\ref{fig:diag1D}(c3) for an example of such locked dynamics. 

Figure~\ref{fig:2D_A2d7}(b) also shows curves S bounding the resonance tongues of some lower-order resonances; we remark that the pairs of curves S bounding narrow resonance tongues are sometimes so close to each other that they appear as single curves in the parameter plane. We find that these pairs of curves S connect different resonance points on the pair of torus bifurcation curves T emerging from the respective Hopf-Hopf bifurcation point HH; specifically a given $p\!\!:\!\!q$ resonance point on one of the curve T connects to a $p\!\!:\!\!(p+q)$ resonance point of the other curve T. To the best of our knowledge, such connecting resonance tongues have been observed only once before, namely in a constructed system of DDEs with two state-dependent delays \cite{CallejaSIADS17}. The Yamada model with a single feedback term and with a constant (non-state-dependent) delay, hence, constitutes the first physically relevant and yet also simplest mathematical model featuring connecting resonance tongues near a Hopf-Hopf bifurcation point. Our results and those in Ref.~[\cite{CallejaSIADS17}] suggest that this phenomenon is generic and related to the occurence of chaotic dynamics near Hopf-Hopf bifurcation. The argument is as follows. Close to the resonance points, the dynamics is either $p\!\!:\!\!q$ or $p\!\!:\!\!(p+q)$ locked on the smooth invariant torus that emerges from the respective torus bifurcation. The fact that the bounding saddle-node bifurcation curves S connect means that the $p\!\!:\!\!q$ locked periodic solutions smoothly evolves into $p\!\!:\!\!(p+q)$ locked solutions. This is perfectly possible in a phase space of dimension at least three, but it cannot happen while these periodic orbits lie on a smooth invariant two-dimensional torus throughout (there cannot be two periodic orbits of different winding number on one and the same two-dimensional torus). Therefore, the smooth invariant torus bifurcating from the respective torus bifurcation T necessarily breaks up at some point when going through the connecting resonance tongue. Such torus break-up can typically lead to the creation of chaotic dynamics through a quasiperiodic route to chaos \cite{Aronson82,green2002bistability}.

The simulations shown in Figure~\ref{fig:diag1D} and, in particular, the appearance of chaotic dynamics can, therefore, be interpreted as follows. For a fixed value of $\tau$ and for small $\kappa$, a stable periodic solution (corresponding to a fixed number of equidistant pulses in the feedback loop) is observed. When $\kappa$ is increased, it looses stability through a torus bifurcation, which leads to the emergence of a stable invariant torus with dynamics that is quasiperiodic or $p\!\!:\!\!q$ locked with very high $q$; see Figure~\ref{fig:diag1D}(c2). As $\kappa$ is increased, the rotation number on the torus evolves and different resonance tongues are crossed, only the larger of which with low $q$ can be identified in Figure~\ref{fig:diag1D}(a); see Figure~\ref{fig:diag1D}(c3) for an example of such locked dynamics. As $\kappa$ is increased further, the torus breaks up and chaotic dynamics is observed over a large range of $\kappa$; as is illustrated in Figure~\ref{fig:diag1D}(c4,c5). Indeed, this overall transition provides an explanation for experimentally observed chaotic pulsing regimes as in Figure~\ref{fig:expe}(d), which are reported here in this microlaser system for the first time.

\subsection{Influence of decreasing the pump parameter $A$}
\label{sec:pump}

We now investigate how the structure of resonance tongues highlighted in Figure~\ref{fig:2D_A2d7}(b) evolves when the pump parameter $A$ is decreased down towards a value $A=2$. The motivation for this is two-fold. First, in the experiment, the pump parameter is related to the amount of energy provided to the microlaser through optical pumping and, as such, is the main control parameter in the experiment. Second, a recent study has shown that for $A=2$ multistable periodic regimes are found, which include non-equidistant pulses in the feedback cavity \cite{TerrienPRE21}. Each of these periodic non-equidistant pulsing regimes exists in a very large region of the $(\tau,\kappa)$-plane of feedback parameters, which is bounded by saddle-node bifurcation curves S of periodic solutions. It has been suggested that these regions are in fact very large resonance tongues, and we now investigate the mechanism for the emergence of such unusally large locking regions. For clarity, we focus on resonance tongues associated with one pair of torus bifurcation curves only, namely the one that emerges from the Hopf-Hopf point HH with the smallest value of $\tau$ in Figure~\ref{fig:2D_A2d7}(b).

\subsubsection{Rotation number near Hopf-Hopf point}
\label{sec:resHH}

\begin{figure}[t!]
\includegraphics[width=\linewidth]{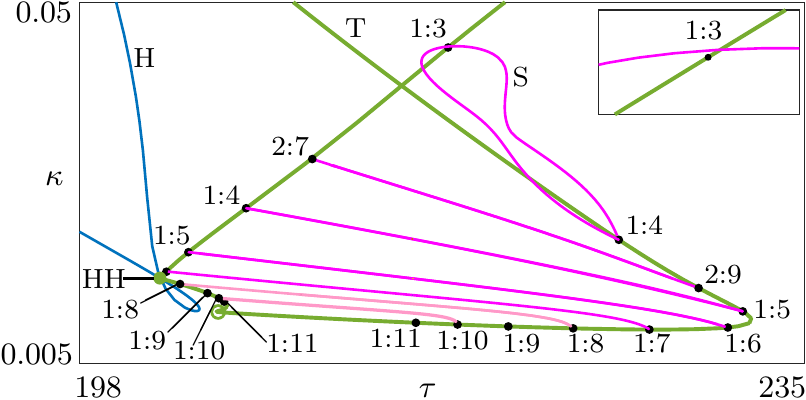}
\caption{Bifurcation diagram in the $(\tau,\kappa)$-plane for $A=2.42$ near the leftmost Hopf-Hopf bifurcation point HH (green dot), showing the self-intersecting Hopf bifurcation curve H (blue), the emerging pair of torus bifurcation curves T (green), and curves S (magenta) of saddle-node bifurcations of periodic orbits that bound resonance tongues arising from selected resonance points (black dots); the open circle on the lower curve T indicates the minimum of the rotation number. Curves S near this minimum that bound resonance tongues that connect two $p\!\!:\!\!q$ resonance points are shown in lighter magenta. The inset is an enlargement near the $1\!\!:\!\!3$ point.}
\label{fig:2D_A2d42}
\end{figure}

Figure~\ref{fig:2D_A2d42} shows the bifurcation diagram of \eqref{eq:yam} for $A=2.42$ in the relevant region of the $(\tau,\kappa)$-plane, featuring the Hopf bifurcation curve H with the Hopf-Hopf bifurcation point HH, the associated two torus bifurcation curves T, and resonance tongues arising from selected $p\!\!:\!\!q$ resonance points; the latter are bounded by pairs of curves S of saddle-node bifurcations of periodic orbits. The lower curve T in Figure~\ref{fig:2D_A2d42} has two fold points (with respect to $\tau$) close to the Hopf-Hopf point HH; moreover, the rotation number $\alpha$ first decreases along this curve away from HH, reaches a minimum of below $\frac{1}{11}$ near (but not exactly at) the leftmost sharp fold, and subsequently increases again. This is in contrast to the properties of the corresponding lower torus bifurcation curves T for $A=2.7$ (see the leftmost part of Figure~\ref{fig:2D_A2d7}(b)), which lacks the leftmost sharp fold and along which  $\alpha$ increases monotonically from the point HH. The non-monotonicity of $\alpha$ along this curve T in Figure~\ref{fig:2D_A2d42} is associated with an intriguing phenomenon: the lighter colored saddle-node curves S connect $p\!\!:\!\!q$ resonance points on the decreasing part to the left of the minimum with $p\!\!:\!\!q$ resonance points on the increasing part to the right of the minimum. Beyond HH, on the other hand, we find connections by pairs of curves S between $p\!\!:\!\!q$ on the upper curve T and $p\!\!:\!\!(p+q)$ resonance points on the other curve T; this situation is as discussed above for $A=2.7$. Another interesting difference with the case $A=2.7$ lies in the shape of the $1\!\!:\!\!4$ resonance tongue in Figure~\ref{fig:2D_A2d42}. Notice that the $1\!\!:\!\!4$ locking region in the $(\tau,\kappa)$-plane is significantly larger than for the other resonances, and it is bounded by a single curve S; as the inset in Figure~\ref{fig:2D_A2d42} shows, this curve passes very close to the $1\!\!:\!\!3$ resonance point, but does not connect to this point.

\subsubsection{Disappearance of Hopf-Hopf point}
\label{sec:disappearHH}

\begin{figure}
\includegraphics[width=\linewidth]{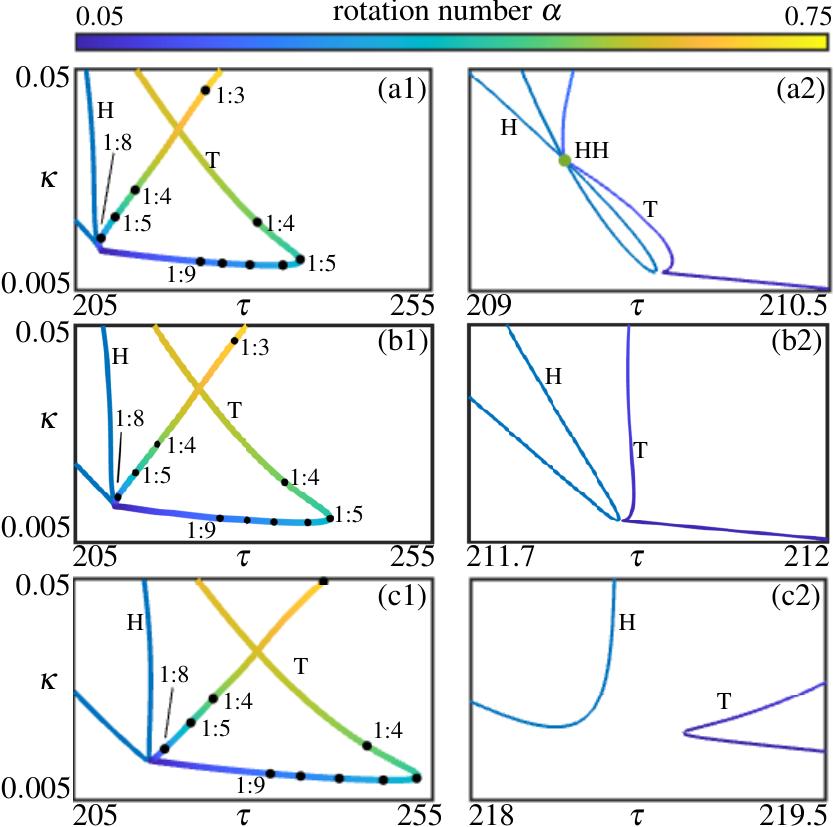}
\caption{Bifurcation diagram in the ($\tau$,$\kappa$)-plane illustrating the disappearance of the Hopf-Hopf bifurcation point HH from Figure~\ref{fig:2D_A2d42} for $A=2.39$ (a), $A=2.38$ (b) and $A=2.35$ (c). Panels (a1)--(c1) show a larger view of the Hopf bifurcation curve H and the torus bifurcation curves T, where the rotation number $\alpha$ along T is indicated according to the color bar and low-order resonance points are marked. Panels (a2)--(c2) are further enlargements.}
\label{fig:torus}
\end{figure}

The Hopf-Hopf bifurcation point HH in Figure~\ref{fig:2D_A2d42} arises because the curve H has a little loop. Comparison with Figure~\ref{fig:2D_A2d7}(b) shows that this loop becomes smaller when the pump parameter $A$ is decreased. Figure~\ref{fig:torus} shows that the loop and the point HH disappear when $A$ is decreased further; moreover, it explains what this means for the associated curves T of torus bifurcation. The rotation number $\alpha$ along the torus bifurcation curves is represented by a color scale in the larger view of panels (a1) to (c1), while panels (a2) to (c2) are further enlargements. The observed configuration in Figure~\ref{fig:torus}(a) for $A=2.39$ is topologically as that for $A=2.42$ in Figure~\ref{fig:2D_A2d42}. However, the point HH in Figure~\ref{fig:torus}(a) is now at a lower value of $\kappa$ and the loop of the curve H is even smaller; similarly, the part of the lower torus curve T along which the rotation number $\alpha$ decreases is significantly smaller than for $A=2.42$. Figure~\ref{fig:torus}(b) for $A=2.38$ is just past the transition of codimension three where the loop disappears. It suggests that this type of degenerate Hopf-Hopf bifurcation occurs when the point HH and the leftmost fold of the lower curve T coincide at a cusp point of the curve H. As a result, for lower values of $A$, as can be observed clearly in Figure~\ref{fig:torus}(c) for $A=2.35$, the two torus bifurcation curves have merged into a single smooth curve T, which is now not connected to the Hopf bifurcation curve. Notice that the rotation number $\alpha$ changes smoothly along this single torus bifurcation curve and has a minimum close to the leftmost fold point of T.

\subsubsection{Growth and merging of resonance regions}
\label{sec:resgrowth}

\begin{figure}[t!]
\includegraphics[width=\linewidth]{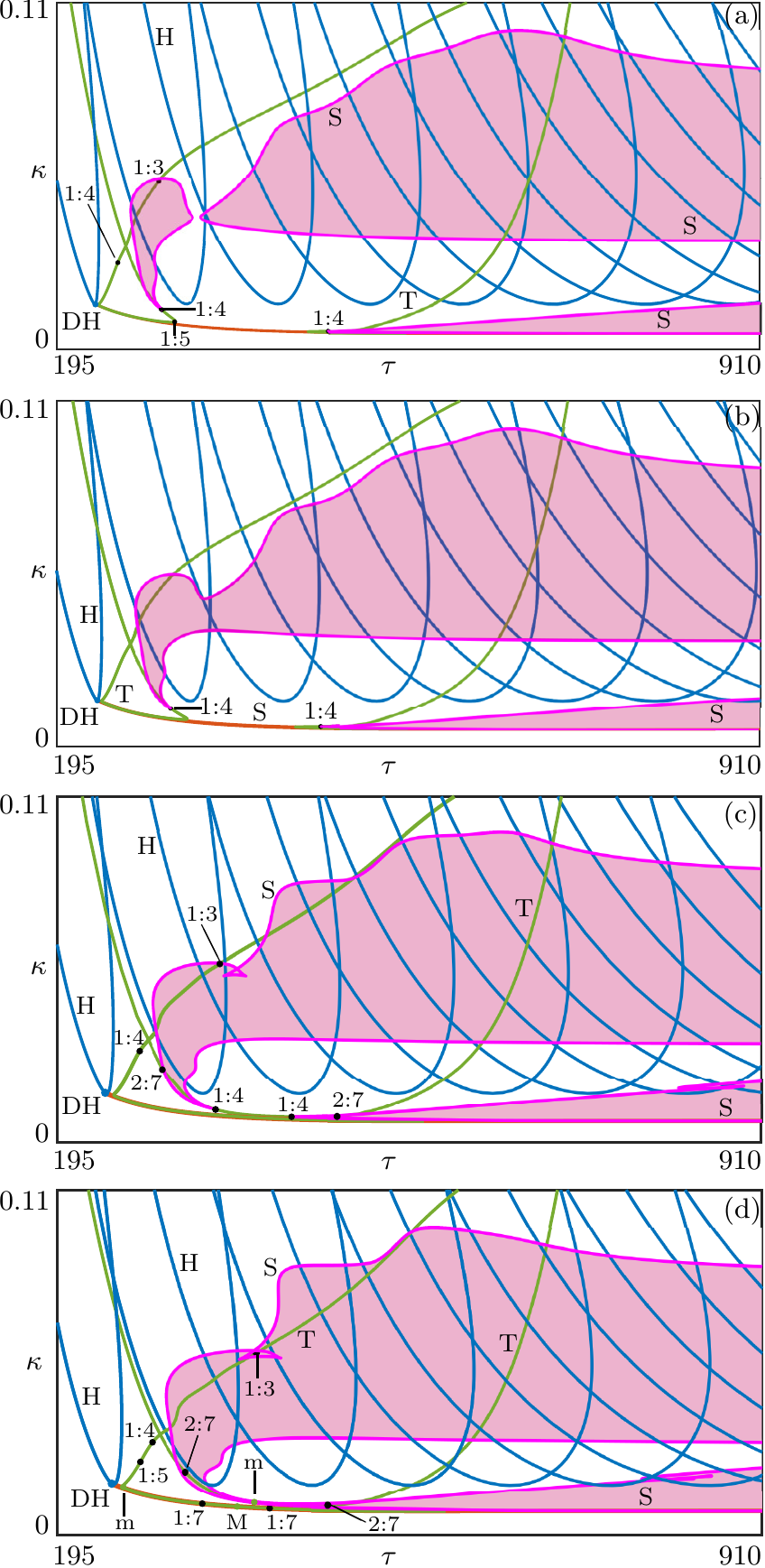}
\caption{Bifurcation diagram of \eqref{eq:yam} in the $(\tau,\kappa)$-plane for $A=2.29$ (a), $A=2.28$ (b), $A= 2.25$ (c), and $A=2.23$ (d), where the region of stable four-pulse solutions is shaded; compare with Figure~\ref{fig:2D_A2d42}.} 
\label{fig:allA}
\end{figure}

We now show how the resonance tongues identified thus far grow into large and experimentally relevant regions of different observable locked periodic dynamics. To this end, Figure~\ref{fig:allA} shows how the bifurcation diagram of \eqref{eq:yam} in the relevant larger region of the $(\tau,\kappa)$-plane changes when the pump parameter $A$ is decreased from $A=2.29$ in panel (a) to $A=2.23$ in panel (d). The focus here is on curves S of saddle-node bifurcation of periodic orbits that bound the (shaded) regions with stable $1\!\!:\!\!4$ dynamics, that is, where one observes periodic pulsing with four (non-equidistant) pulses in the feedback cavity. 

Figure~\ref{fig:allA}(a) for $A=2.29$ shows that the $1\!\!:\!\!4$ resonance tongue bounded by a single curve S, which was identified in Figure~\ref{fig:2D_A2d42}, has grown; note that this curve S in Figure~\ref{fig:allA}(a) still connects to and from a $1\!\!:\!\!4$ resonance point on the torus bifurcation curve T associated with the leftmost point HH that just disappeared. Close by, there is another, large $1\!\!:\!\!4$ resonance region that is bounded by another curve S and exists for larger values of $\tau$; we remark that this resonance region does not connect to any resonance point and extends to very large values of $\tau$ (as we checked by computing its bounding curve S beyond the shown range). 
When the pump parameter is changed to $A=2.28$ as shown in Figure~\ref{fig:allA}(b), these two resonance regions have merged into a single and very large $1\!\!:\!\!4$ locking region emerging from the left-most $1\!\!:\!\!4$ resonance point, which happens to be near the right-most fold on the associated curve T. This change occurs in a saddle transition of the bounding curves S of saddle-node bifurcations of periodic orbits, where the two nearby fold points of the two curves S meet and the curves connect differently. 

Panels (a) and (b) of Figure~\ref{fig:allA} also show a further locking region that emerges from a $1\!\!:\!\!4$ resonance point on a second torus bifurcation curve T, which has a fold with respect to $\tau$ and then terminates for low values of $\kappa$ at a $1\!\!:\!\!1$ resonance point on the shown curve S of saddle-node bifurcation of the (equidistant) periodic pulsing solutions \cite{TerrienSIADS17}; note that this curve S emerges near the leftmost minimum of the curve H from a degenerate Hopf bifurcation point DH. The two fold points of the two curves T move closer to each other with decreasing $A$, and in Figure~\ref{fig:allA}(c) for $A=2.25$ they now connect differently; also after a saddle-transition, now of the curves T. This means, in particular, that the two $1\!\!:\!\!4$ resonance points now lie on the same curve T. As $A$ is decreased further, these two points meet and then disappear, as is shown in Figure~\ref{fig:allA}(d). Simultaneousy, the two resonance tongues meet, and their bounding curves S connect differently. This leads to a very large resonance region (shaded) of $1\!\!:\!\!4$ locked pulsing, bounded by a curve S, which is not connected anymore to any of the torus bifurcation curves found in the considered parameter range. Inside this large four-pulse region, one finds a smaller region, bounded by a second curve S, where the locked solution does not exist. Notice further that pairs of cusp points have occured in Figure~\ref{fig:allA}(c) and (d) on the upper curve S and in panel (d) on the lower intermediate curve S; these cusp points appear locally in codimension-three swallowtail bifurcations and lead to small regions with additional $1\!\!:\!\!4$ locked solutions that are not discussed further in this article. 

\begin{figure}[t!]
\includegraphics[width=\linewidth]{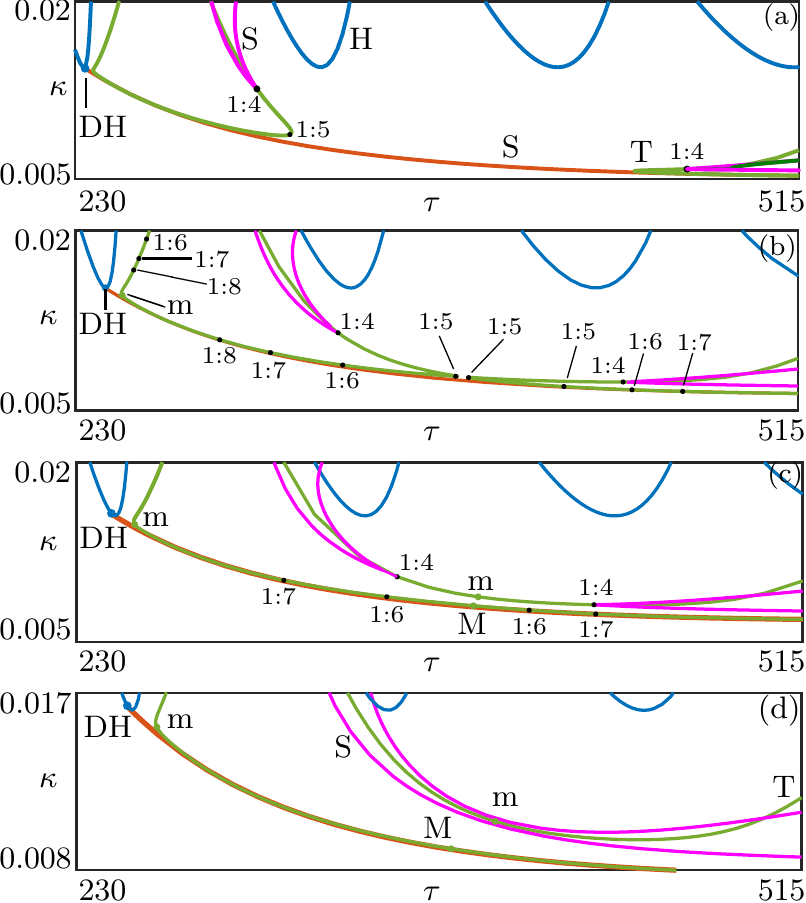}
\caption{Disconnecting $1\!\!:\!\!4$ resonance tongues of \eqref{eq:yam} in the relevant region of the $(\tau,\kappa)$-plane, showing the relevant bifurcation curves for $A=2.29$ (a), $A=2.2605$ (b), $A= 2.25$ (c), and $A=2.23$ (d); compare with Figure~\ref{fig:allA}.} 
\label{fig:allAzoom}
\end{figure}

The changes to the curves T and the bounding curves S of the $1\!\!:\!\!4$ resonance regions occur over a small parameter range of $A$ and are difficult to see on the scale of Figure~\ref{fig:allA}. This is why they are illustrated further in Figure~\ref{fig:allAzoom} with enlargements near the two $1\!\!:\!\!4$ resonance points involved. Figure~\ref{fig:allAzoom}(a) for $A=2.29$ (an enlargement of Figure~\ref{fig:allA}(a)) clearly shows the respective two resonance tongues emerging from the $1\!\!:\!\!4$ resonance points near the folds of the two different torus bifurcation curves T. Panel (b) for $A=2.2605$ is just before the saddle transition of these curves T: the two folds are now extremely close together and there are two $1\!\!:\!\!5$ resonance points practically on them. When the curves T are connected differently, just after the saddle transition for $A= 2.25$ as shown in Figure~\ref{fig:allAzoom}(c) (an enlargement of Figure~\ref{fig:allA}(c)), the two $1\!\!:\!\!4$ resonance points are now clearly on the same curve T. This means, in particular, that there must be an extremum of the rotation number. Different resonance points are indicated along the curves T to give information on the rotation number along them, and this shows that the point m is a minimum near $\frac{1}{5}$. Similarly, there is a maximum $M$ of the rotation number, also near $\frac{1}{5}$, on the other torus bifurcation curve T in panel (c). The value of the rotation number at the minimum m increases as $A$ is decreased; in the process, the two $1\!\!:\!\!4$ resonance points move towards each other until they collide and disappear when the minimum moves through $\frac{1}{4}$, as discussed previously. The result is the situation for $A=2.23$ in Figure~\ref{fig:allAzoom}(d) (an enlargement of Figure~\ref{fig:allA}(d)), where the bounding curves S of the associated resonance tongues are now connected differently: they lie either side of the point m on T and no longer interact with this curve of torus bifurcation. Note that, throughout in Figure~\ref{fig:allAzoom}, the lower parts of curves T are very close to the curve S that emerges from the point DH near the minimum of the Hopf bifurcation on the left.

Finally, as $A$ is decreased further, the region without locking inside the large $1\!\!:\!\!4$ pulsing region in Figure~\ref{fig:allA} shrinks and eventually disappears from the shown range of $\tau$ through a minimax transition of the bounding curve S (this is not illustrated here). Overall, this leads to a single, very large locking region, which is not connected to any resonance point along a torus bifurcation curve and which extends to very large values of the delay time $\tau$. Similar transitions of other resonance tongues exist as well (but are not shown in Figures~\ref{fig:allA} and~\ref{fig:allAzoom}); this agrees with the observation of associated regions with stable periodic solutions with different numbers of non-equidistant pulses \cite{TerrienPRE21}. From a practical point of view, this means that, for these lower values of the pump parameter $A$, locked periodic solutions corresponding to such non-equidistant pulses in the external feedback cavity can be observed over very large ranges in the $(\tau,\kappa)$-plane. This agrees with the behaviour observed recently in both the model and an actual experiment \cite{TerrienPRE21}, which confirms that the non-equidistant pulsing periodic regimes originate in a resonance phenomenon. 

While the resulting regions of stable pulsing solutions are large in the $(\tau,\kappa)$-plane, we stress that the transitions of resonance tongues that generate them happen in a very narrow range of $A$. Hence, the system is highly sensitive to small changes of the pump parameter. This is of practical importance since the pump parameter $A$ is a main control parameter in the actual experiment: even a small change in $A$ may result in the disappearance of the locked solution (corresponding to non-equidistant pulses in the feedback cavity) and, as a consequence, the appearance of qualitatively different long-term dynamics\cite{TerrienPRR20}. This observed change with $A$ appears to be sudden from an experimental and practical perspective; mathematically, on the other hand, it is explained by the sequence of transitions of resonance tongues presented above.

\subsection{Disconnecting and disappearing resonance tongues near extrema of the rotation number}
\label{sec:merge_dis}

\begin{figure}
\includegraphics[width=\linewidth]{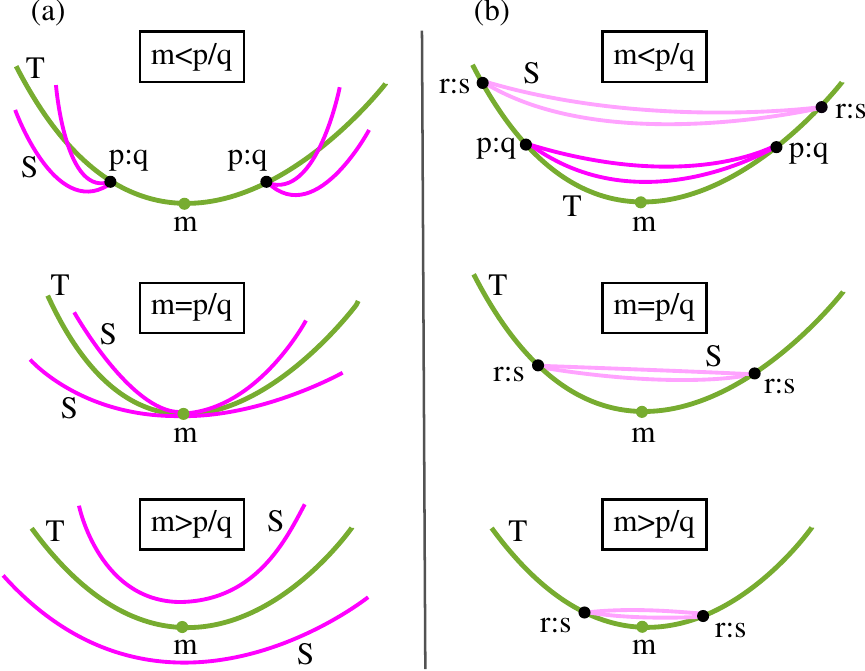}
\caption{Sketches of the two scenarios for disconnecting (a) and disappearing (b) $p\!\!:\!\!q$ resonance tongues near an extremum of the rotation number on a curve T of torus bifurcation. Shown is the case for a minumum m where, from top to bottom, ${\rm m} < \frac{p}{q}$, ${\rm m} = \frac{p}{q}$ and ${\rm m} > \frac{p}{q}$; column (b) also shows a connecting $r\!:\!s$ resonance tongue with ${\rm m} < \frac{r}{s}$.}
\label{fig:sketch_merg_dis}
\end{figure}

While Figures~\ref{fig:allA} and~\ref{fig:allAzoom} illustrate the case of $1\!\!:\!\!4$ resonance, disconnecting resonance tongues are a generic phenomenon that occurs for any pair of $p\!\!:\!\!q$ resonance points near an extremum of the rotation number on a torus bifurcation curve T. We now discuss in more detail the disappearance of pairs of resonance points and the consequences for the associated resonance tongues. To our knowledge, this generic phenomenon of codimension three (meaning that it happens at a specific point in a three-dimensional parameter space, as the $(\tau,\kappa,A)$-space of \eqref{eq:yam} considered here) has not been reported in previous literature. In fact, there are actually two different cases, as is sketched in Figure~\ref{fig:sketch_merg_dis} for the case that the extremum on T is a minimum m in a two-parameter plane. In both columns~(a) and~(b) of Figure~\ref{fig:sketch_merg_dis} a pair of $p\!\!:\!\!q$ resonance points with ${\rm m} < \frac{p}{q}$ collides when $m$ increases and reaches the value ${\rm m} = \frac{p}{q}$ as a third parameter is changed, and subsequently does not exist any longer on the curve T when ${\rm m} > \frac{p}{q}$. As a function of the third parameter ($A$ for \eqref{eq:yam}), this corresponds to a fold point of $p\!\!:\!\!q$ resonance points on the curve T. Column~(a) shows the situation encountered in Section~\ref{sec:resgrowth}, where the resonance tongues `point away' from the minimum m; as a result, the pairs of bounding curves S then meet at m when ${\rm m} = \frac{p}{q}$ and form two curves on either side of the curve T for ${\rm m} > \frac{p}{q}$. Since the bounding curves S are then no longer connected with T at a point of resonance, we refer to this generic case as that of \emph{disconnecting resonance tongues}. In contrast, the $p\!\!:\!\!q$ resonance tongues in Figure~\ref{fig:sketch_merg_dis}(b) connect near the minimum m, meaning that the curves S bound a single small resonance region. This $p\!\!:\!\!q$ resonance region shrinks to a point when ${\rm m} = \frac{p}{q}$ and has disappeared for ${\rm m} > \frac{p}{q}$; left are then only $r\!:\!s$ resonance tongue with ${\rm m} < \frac{r}{s}$, of which one is shown in the sketch. We refer to this generic case as that of \emph{disappearing resonance tongues}.

Figure~\ref{fig:sketch_merg_dis} illustrates the codimension-three transitions of disconnecting and disappearing resonance tongues locally for a single pair of resonance points. As the value of m increases smoothly and monotonically with a third parameter, the same local scenario takes place at any point on the torus bifurcation T where the rotation number is rational. Since this set is dense in T, any change of the parameter leads to infinitely many disconnecting or disappearing resonance tongues. In other words, Figure~\ref{fig:sketch_merg_dis} represents a countably infinite sequence of such transitions of resonance tongues as the third parameter is changed and the minimum m increases. We remark that the case of a maximum M on T is completely analogous: here the resonance tongues merge or disappear at rational points when M decreases with a third parameter. With this interpretation, we claim that columns~(a) and~(b) of Figure~\ref{fig:sketch_merg_dis} represent the unfolding of the two cases of disconnecting and of disappearing resonance tongues. 

\begin{figure}
\includegraphics[width=\linewidth]{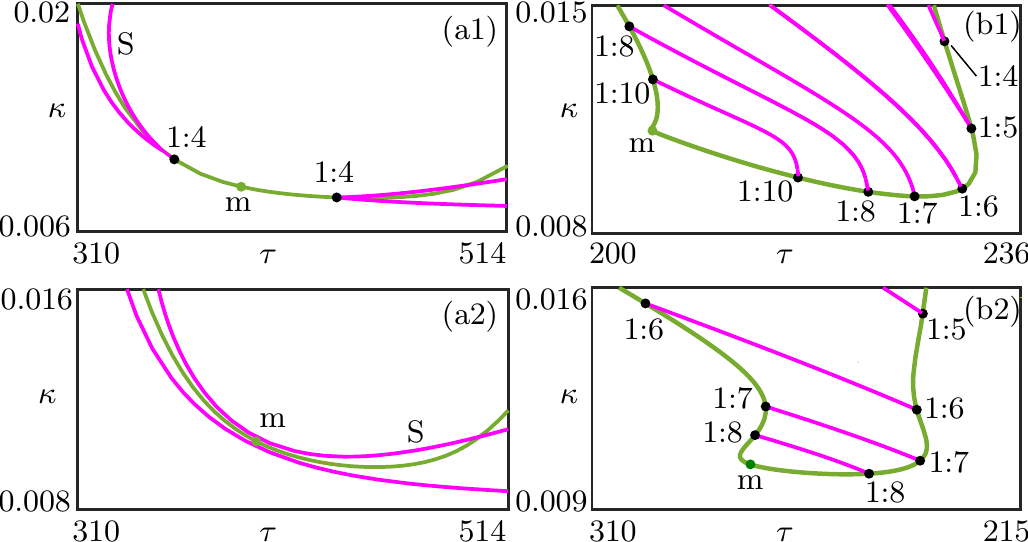}
\caption{Disconnecting (a) and disappearing (b) resonance tongues in the $(\tau,\kappa)$-plane of \eqref{eq:yam} near a minimum m of the rotation number on the torus bifurcation T; here $A = 2.25$ with ${\rm m} = 0.234$ (a1), $A = 2.23$ with ${\rm m} = 0.257$ (a2), $A = 2.42$ with ${\rm m} = 0.032$ (b1), and $A = 2.46$ with ${\rm m} = 0.104$ (b2). Compare panels (a) with Figure~\ref{fig:allAzoom}(c)--(d) and panels (b) with Figure~\ref{fig:2D_A2d42}.}
\label{fig:data_merg_dis}
\end{figure}

Proving that these unfoldings are generic is a challenge that is left for future work, especially since it would require the simpler setting of an ODE rather than a DDE as studied here. However, the correctness of these two unfoldings is supported by numerical evidence in the $(\tau,\kappa)$-plane of \eqref{eq:yam}. Specifically, Figure~\ref{fig:data_merg_dis} highlights the local situations before and after the respective transition for the disconnecting resonance tongues in column~(a) and for disappearing resonance tongues in column~(b);  here the parameter $A$ changes the value of the minimum m on the curve T as indicated. Figure~\ref{fig:data_merg_dis}(a1) and~(a2) show a further enlargement that illustrates how the two shown $1\!\!:\!\!4$ resonance tongues from Figure~\ref{fig:allAzoom} disconnect from T; compare with Figure~\ref{fig:sketch_merg_dis}(a). Similarly, panels~(b1) and~(b2) show the disappearance of a resonance tongue that connects two $1\!\!:\!\!10$ resonance points on the curve T; notice that the other shown resonance points and tongues have moved closer to m and will disappear when $A$ is increased further. Note that these resonance tongues disappear at the minimum near the Hopf-Hopf bifurcation point HH shown in Figure~\ref{fig:2D_A2d42}.

\begin{figure}
\includegraphics[width=\linewidth]{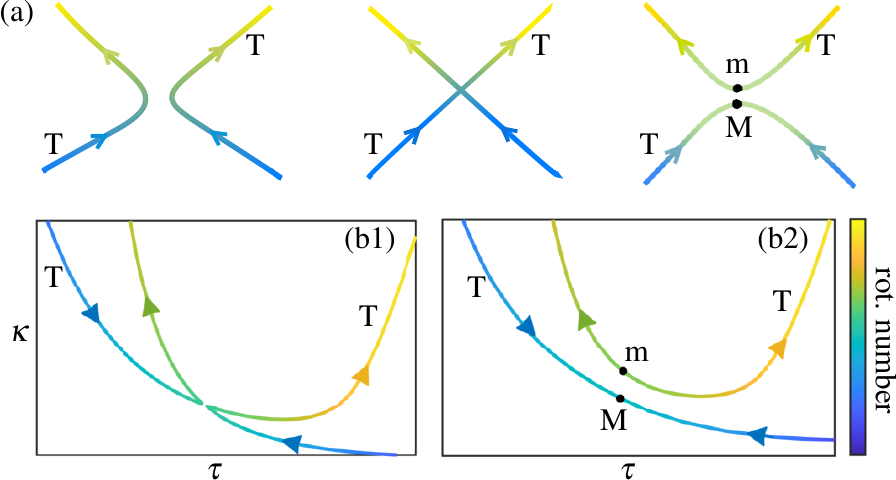}
\caption{Saddle transition of torus bifurcation curves generating a minimum m and a maximum M of the rotation number. Panels (a) are sketches, while panels (b) show the $(\tau,\kappa)$-plane of \eqref{eq:yam} just before (b1) and after (b2) the saddle transition, for $A = 2.25$ and $A = 2.23$, respectively. Compare with Figure~\ref{fig:allAzoom}(b)--(c).}
\label{fig:torus_minmax}
\end{figure}

To complete this section, we present how extrema arise on a torus bifurcation curve in a two-parameter plane. As row~(a) of Figure~\ref{fig:torus_minmax} shows, the generic mechanism is that of a saddle-transition of a pair of torus bifurcation curves, which necessarily needs to respect and agree with the respective changes of the rotation number along these curves. Note that generically the rotation number is not constant, and Figure~\ref{fig:torus_minmax}(a) shows the case that it increases as indicated by the arrows on the curves labeled T. As a third parameter is changed, the saddle transition occurs: here the two curves meet as shown at a singular point (which is a saddle point of the surface of torus bifurcations in the three-dimensional parameter space). The locus of torus bifurcations connects differently past the saddle transition, as is shown. On the level of just the curves T, this is a standard way of reconnecting the respective branches of torus bifurcation. However, we can conclude more here: since the rotation number is unique and common to all branches at the singular point, a minimum m on one curve T and a maximum M on the other curve T are necessarily created in the process. The sketch in Figure~\ref{fig:torus_minmax}(a) represents the generic case of the saddle transition of a locus of torus bifurcations. As we discussed in Section~\ref{sec:resgrowth}, such a saddle transition occurs in the $(\tau,\kappa)$-plane of \eqref{eq:yam}, namely between the two situations shown in Figure~\ref{fig:allAzoom}(b) and~(c). This is illustrated more clearly in Figure~\ref{fig:torus_minmax}(b) with enlargements near the saddle transition that also shows the rotation number along the respective curves T. Notice that the situation in panel~(b1) is very close to the saddle-transitions, where the curves T connect.

\section{Conclusion}
\label{sec:concl}

We investigated the emergence of complex multi-frequency dynamics in an excitable system subject a delayed feedback loop --- specifically a microlaser with optical feedback from an external mirror. The Yamada model with delayed feedback has been shown to reproduce accurately a range of phenomena observed experimentally in an excitable microlaser subject to delayed optical feedback. This includes quasiperiodic dynamics, locked periodic dynamics on a torus and chaotic dynamics. A bifurcation analysis in the physically relevant variables, the feedback delay $\tau$, the feedback strength $\kappa$ and the pump parameter $A$, unveiled the key role played by resonance tongues in the $(\tau,\kappa)$-plane, which arise from torus bifurcation curves that, in turn, emerge from Hopf-Hopf bifurcation points on the main curve of Hopf bifurcations. In particular, the emergence of chaotic dynamics in the DDE results from the mechanism of torus break-up. This is associated with the observation of \textit{connecting} resonance tongues, which connect $p\!\!:\!\!q$ resonance points with $p\!\!:\!\!(p+q)$resonance points. 

Our results confirm that non-equidistant pulsing periodic solutions originate in resonance phenomena. Unexpectedly large locking regions in which such pulsing regimes are observed emerge in the $(\tau,\kappa)$-plane from initially much smaller resonance tongues in an intriguing sequence of transitions when the pump parameter $A$ is changed. In particular, this involves several re-arrangements of pairs of bifurcation curves in saddle transitions, as well as pairs of $p\!\!:\!\!q$ resonance points coming together at points on a torus bifurcation curve where the rotation number has an extremum. We presented (conjectural) unfoldings of the two generic cases of disconneting and of disappearing resonance tongues, and also showed how extrema on torus bifurcation curves arise. In the Yamada model with delayed feedback the respective transitions all occur in quick succession, showing that the system is very sensitive to small changes in the pump parameter $A$. Because $A$ is a main control parameter this is of practical interest in an actual experiment --- the bifurcation analysis presented here explains the consistency and genericity of what may otherwise be interpreted as a non-obvious sudden jump from stable single-pulse behavior to pulsing with several equidistant and/or non-equidistant pulses.

From a more general point of view, the DDE model considered here has only two main ingredients: excitability and feedback. As such, our results are expected to be of more general practical interest, beyond the particular laser device considered in this paper. Examples with these same ingredients are other optical systems with different types of feedback, as well as biological systems with delayed feedback, which show considerable similarities, as confirmed by a recent study with an excitable cell \cite{WedgwoodJRS21}.

\section*{Acknowledgments}
VAP, RB, IS, GB, KP and SB acknowledge partial support from the French Network of Nanotechology facilities Renatech.



\begin{thebibliography}{10}

\bibitem{ArecchiPRA92}
FT~Arecchi, G~Giacomelli, A~Lapucci, and R~Meucci.
\newblock Two-dimensional representation of a delayed dynamical system.
\newblock {\em Physical Review A}, 45(7):R4225, 1992.

\bibitem{Aronson82}
DG~Aronson, MA~Chory, GR~Hall, and Richard~P McGehee.
\newblock Bifurcations from an invariant circle for two-parameter families of
  maps of the plane: a computer-assisted study.
\newblock {\em Communications in Mathematical Physics}, 83(3):303--354, 1982.

\bibitem{BarbayOL11}
S~Barbay, R~Kuszelewicz, and AM~Yacomotti.
\newblock Excitability in a semiconductor laser with saturable absorber.
\newblock {\em Opt. Lett.}, 36(23):4476--4478, Dec 2011.

\bibitem{CallejaSIADS17}
RC~Calleja, AR~Humphries, and B~Krauskopf.
\newblock Resonance phenomena in a scalar delay differential equation with two
  state-dependent delays.
\newblock {\em SIAM Journal on Applied Dynamical Systems}, 16(3):1474--1513,
  2017.

\bibitem{DubbeldamOC99}
JLA~Dubbeldam and B~Krauskopf.
\newblock Self-pulsations of lasers with saturable absorber : dynamics and
  bifurcations.
\newblock {\em Opt. Commun.}, 159:325, 1999.

\bibitem{ElsassEPJD10}
T~Elsass, K~Gauthron, G~Beaudoin, I~Sagnes, R~Kuszelewicz, and S.~Barbay.
\newblock Control of cavity solitons and dynamical states in a monolithic
  vertical cavity laser with saturable absorber.
\newblock {\em Eur. Phys. J. D}, 59, 2010.

\bibitem{Engelborghs00report}
K~Engelborghs, T~Luzyanina, G~Samaey, et~al.
\newblock Dde-biftool: a matlab package for bifurcation analysis of delay
  differential equations.
\newblock {\em TW Report}, 305:1--36, 2000.

\bibitem{EngelborghsACM02}
K~Engelborghs, T~Luzyanina, and D~Roose.
\newblock Numerical bifurcation analysis of delay differential equations using
  dde-biftool.
\newblock {\em ACM Transactions on Mathematical Software (TOMS)}, 28(1):1--21,
  2002.

\bibitem{ErneuxPRE18}
T~Erneux and S~Barbay.
\newblock Two distinct excitable responses for a laser with a saturable
  absorber.
\newblock {\em Phys. Rev. E}, 97(6), jun 2018.

\bibitem{garbin2017interactions}
B~Garbin, J~Javaloyes, S~Barland, and G~Tissoni.
\newblock Interactions and collisions of topological solitons in a
  semiconductor laser with optical injection and feedback.
\newblock {\em Chaos: An Interdisciplinary Journal of Nonlinear Science},
  27(11):114308, 2017.

\bibitem{GarbinNC15}
B~Garbin, J~Javaloyes, G~Tissoni, and S~Barland.
\newblock Topological solitons as addressable phase bits in a driven laser.
\newblock {\em Nat. Commun.}, 6:--, jan 2015.

\bibitem{green2002bistability}
K~Green, B~Krauskopf, and K~Engelborghs.
\newblock Bistability and torus break-up in a semiconductor laser with
  phase-conjugate feedback.
\newblock {\em Physica D: Nonlinear Phenomena}, 173(1-2):114--129, 2002.

\bibitem{IzhikevichBook}
EM~Izhikevich.
\newblock {\em Dynamical Systems in Neuroscience: The Geometry of Excitability
  and Bursting.}
\newblock The MIT press, 2007.

\bibitem{izhikevich2000neural}
EM~Izhikevich.
\newblock Neural excitability, spiking and bursting.
\newblock {\em International journal of bifurcation and chaos},
  10(06):1171--1266, 2000.

\bibitem{keane2015delayed}
A~Keane, B~Krauskopf, and C~Postlethwaite.
\newblock Delayed feedback versus seasonal forcing: resonance phenomena in an
  el ni{\~n}o southern oscillation model.
\newblock {\em SIAM Journal on Applied Dynamical Systems}, 14(3):1229--1257,
  2015.

\bibitem{KrauskopfWalker}
B~Krauskopf and JJ~Walker.
\newblock {\em Bifurcation Study of a Semiconductor Laser with Saturable
  Absorber and Delayed Optical Feedback}, pages 161--181.
\newblock Wiley-VCH Verlag GmbH \& Co. KGaA, 2012.

\bibitem{Kuznetsov2013elements}
YA~Kuznetsov.
\newblock {\em Elements of applied bifurcation theory}, volume 112.
\newblock Springer Science \& Business Media, 2013.

\bibitem{munsberg2020topological}
L~Munsberg, J~Javaloyes, and SV~Gurevich.
\newblock Topological localized states in the time delayed adler model:
  Bifurcation analysis and interaction law.
\newblock {\em Chaos: An Interdisciplinary Journal of Nonlinear Science},
  30(6):063137, 2020.

\bibitem{otupiri2020yamada}
R~Otupiri, B~Krauskopf, and NGR~Broderick.
\newblock The yamada model for a self-pulsing laser: Bifurcation structure for
  nonidentical decay times of gain and absorber.
\newblock {\em International Journal of Bifurcation and Chaos}, 30(14):2030039,
  2020.

\bibitem{pammi2019photonic}
VA~Pammi, K~Alfaro-Bittner, MG~Clerc, and S~
  Barbay.
\newblock Photonic computing with single and coupled spiking micropillar
  lasers.
\newblock {\em IEEE Journal of Selected Topics in Quantum Electronics},
  26(1):1--7, 2019.

\bibitem{RomeiraNSR16}
B~Romeira, R~Av\'o, JML~ Figueiredo, S~Barland, and J~Javaloyes.
\newblock Regenerative memory in time-delayed neuromorphic photonic resonators.
\newblock {\em Sci. Rep.}, 6, 2016.

\bibitem{roose2007continuation}
D~Roose and R~ Szalai.
\newblock Continuation and bifurcation analysis of delay differential
  equations.
\newblock In {\em Numerical continuation methods for dynamical systems}, pages
  359--399. Springer, 2007.

\bibitem{ruschel2020limits}
S~Ruschel, B~Krauskopf, and NGR~Broderick.
\newblock The limits of sustained self-excitation and stable periodic pulse
  trains in the yamada model with delayed optical feedback.
\newblock {\em Chaos: An Interdisciplinary Journal of Nonlinear Science},
  30(9):093101, 2020.

\bibitem{seidel2022influence}
Thomas~G Seidel, Julien Javaloyes, and Svetlana~V Gurevich.
\newblock Influence of time-delayed feedback on the dynamics of temporal
  localized structures in passively mode-locked semiconductor lasers.
\newblock {\em Chaos: An Interdisciplinary Journal of Nonlinear Science},
  32(3):033102, 2022.

\bibitem{SelmiPRL14}
F.~Selmi, R.~Braive, G.~Beaudoin, I.~Sagnes, R.~Kuszelewicz, and S.~Barbay.
\newblock Relative refractory period in an excitable semiconductor laser.
\newblock {\em Phys. Rev. Lett.}, 112:183902, May 2014.

\bibitem{SelmiPRE16}
F~Selmi, R~Braive, G~Beaudoin, I.-~Sagnes, R~Kuszelewicz, T~Erneux, and
  S~Barbay.
\newblock Spike latency and response properties of an excitable micropillar
  laser.
\newblock {\em Phys. Rev. E}, 94:042219, Oct 2016.

\bibitem{shastri2021photonics}
BJ~Shastri, AN~Tait, T~Ferreira de Lima, WHP~Pernice, Harish Bhaskaran, CD~Wright, and PR~Prucnal.
\newblock Photonics for artificial intelligence and neuromorphic computing.
\newblock {\em Nature Photonics}, 15(2):102--114, 2021.

\bibitem{sieber2014dde}
J~ Sieber, K~Engelborghs, T~ Luzyanina, G~Samaey, and D~
  Roose.
\newblock Dde-biftool manual-bifurcation analysis of delay differential
  equations.
\newblock {\em arXiv preprint arXiv:1406.7144}, 2014.

\bibitem{TerrienSIADS17}
S~Terrien, B~Krauskopf, and NGR~ Broderick.
\newblock Bifurcation analysis of the yamada model for a pulsing semiconductor
  laser with saturable absorber and delayed optical feedback.
\newblock {\em SIAM J. Appl. Dyn. Sys.}, 16(2):771--801, 2017.

\bibitem{TerrienPRA17}
S~Terrien, B~ Krauskopf, NGR~. Broderick, L~ Andr{\'{e}}oli,
  F~ Selmi, R~ Braive, G~ Beaudoin, I~Sagnes, and
  S~Barbay.
\newblock Asymmetric noise sensitivity of pulse trains in an excitable
  microlaser with delayed optical feedback.
\newblock {\em Phys. Rev. A}, 96(4), oct 2017.

\bibitem{TerrienOL18}
S~ Terrien, B~ Krauskopf, NGR ~Broderick, R~ Braive,
  G~ Beaudoin, I~Sagnes, and S~Barbay.
\newblock Pulse train interaction and control in a microcavity laser with
  delayed optical feedback.
\newblock {\em Opt. Lett.}, 43(13):3013--3016, 2018.

\bibitem{TerrienPRR20}
S~Terrien, VA~ Pammi, NGR~ Broderick, R~ Braive,
  G~Beaudoin, I~ Sagnes, B~ Krauskopf, and S~
  Barbay.
\newblock Equalization of pulse timings in an excitable microlaser system with
  delay.
\newblock {\em Phys. Rev. Research}, 2(2), apr 2020.

\bibitem{TerrienPRE21}
S~Terrien, VA~Pammi, B~Krauskopf, NGR~Broderick, and
  S~Barbay.
\newblock Pulse-timing symmetry breaking in an excitable optical system with
  delay.
\newblock {\em Physical Review E}, 103(1):012210, 2021.

\bibitem{turconi2013control}
M~Turconi, B~Garbin, M~Feyereisen, M~Giudici, and S~Barland.
\newblock Control of excitable pulses in an injection-locked semiconductor
  laser.
\newblock {\em Physical Review E}, 88(2):022923, 2013.

\bibitem{WedgwoodJRS21}
KCA~Wedgwood, P~S{\l}owi{\'n}ski, J~Manson, K~ Tsaneva-Atanasova, and B~Krauskopf.
\newblock Robust spike timing in an excitable cell with delayed feedback.
\newblock {\em Journal of the Royal Society Interface}, 18(177):20210029, 2021.

\bibitem{Yamada93}
M~Yamada.
\newblock A theoretical analysis of self-sustained pulsation phenomena in
  narrow-stripe semiconductor lasers.
\newblock {\em IEEE J. Quantum Electron.}, 29(5):1330--1336, 1993.

\bibitem{yanchuk2009delay}
S~Yanchuk and P~Perlikowski.
\newblock Delay and periodicity.
\newblock {\em Physical Review E}, 79(4):046221, 2009.

\end{thebibliography}
\end{document}